%% Plain TeX
\newcount\secno
\newcount\prmno
\newif\ifnotfound
\newif\iffound

\def\namedef#1{\expandafter\def\csname #1\endcsname}
\def\nameuse#1{\csname #1\endcsname}

\long\def\ifundefined#1#2#3{\expandafter\ifx\csname
  #1\endcsname\relax#2\else#3\fi}
\def\hwrite#1#2{{\let\the=0\edef\next{\write#1{#2}}\next}}

% Working with lists
\toksdef\ta=0 \toksdef\tb=2
\long\def\leftappenditem#1\to#2{\ta={\\{#1}}\tb=\expandafter{#2}%
                                \edef#2{\the\ta\the\tb}}
\long\def\rightappenditem#1\to#2{\ta={\\{#1}}\tb=\expandafter{#2}%
                                \edef#2{\the\tb\the\ta}}

\def\lop#1\to#2{\expandafter\lopoff#1\lopoff#1#2}
\long\def\lopoff\\#1#2\lopoff#3#4{\def#4{#1}\def#3{#2}}

\def\ismember#1\of#2{\foundfalse{\let\given=#1%
    \def\\##1{\def\next{##1}%
    \ifx\next\given{\global\foundtrue}\fi}#2}}

% Les commandes
\def\section#1{\vskip1truecm
               \global\def\currenvir{section}
               \global\advance\secno by1\global\prmno=0
               {\bf \number\secno. {#1}}
               \smallskip}

\def\subsection{\global\def\currenvir{subsection}
                \global\advance\prmno by1
                \ind{ (\number\secno.\number\prmno) }}
\def\subsec{\global\def\currenvir{subsection}
                \global\advance\prmno by1
                { (\number\secno.\number\prmno)\ }}

\def\proclaim#1{\global\advance\prmno by 1
                {\bf #1 \the\secno.\the\prmno$.-$ }}

\long\def\th#1 \enonce#2\endth{%
   \medbreak\proclaim{#1}{\it #2}\global\def\currenvir{th}\smallskip}

\def\rem#1{\global\advance\prmno by 1
{\it #1} \the\secno.\the\prmno$.-$}

% CROSS-REFERENCES

\def\isinlabellist#1\of#2{\notfoundtrue%
   {\def\given{#1}%
    \def\\##1{\def\next{##1}%
    \lop\next\to\za\lop\next\to\zb%
    \ifx\za\given{\zb\global\notfoundfalse}\fi}#2}%
    \ifnotfound{\immediate\write16%
                 {Warning - [Page \the\pageno] {#1} No reference found}}%
                \fi}%
\def\ref#1{\ifx\labellist\empty{\immediate\write16
                 {Warning - No references found at all.}}
               \else{\isinlabellist{#1}\of\labellist}\fi}

\def\newlabel#1#2{\rightappenditem{\\{#1}\\{#2}}\to\labellist}
\def\labellist{}

\def\label#1{%
  \def\given{th}%
  \ifx\given\currenvir%
    {\hwrite\lbl{\string\newlabel{#1}{\number\secno.\number\prmno}}}\fi%
  \def\given{section}%
  \ifx\given\currenvir%
    {\hwrite\lbl{\string\newlabel{#1}{\number\secno}}}\fi%
  \def\given{subsection}%
  \ifx\given\currenvir%
    {\hwrite\lbl{\string\newlabel{#1}{\number\secno.\number\prmno}}}\fi%
  \def\given{subsubsection}%
  \ifx\given\currenvir%
  {\hwrite\lbl{\string%
    \newlabel{#1}{\number\secno.\number\subsecno.\number\subsubsecno}}}\fi
  \ignorespaces}

\newwrite\lbl
\def\begin{\newlabel{mori}{1.2}
\newlabel{dp}{1.3}
\newlabel{bireg}{1.4}
\newlabel{ii}{1.5}
\newlabel{dj2}{1.6}
\newlabel{modele}{2.1}
\newlabel{fixed}{2.2}
\newlabel{equi}{2.3}}

%% Fin de la numerotation automatique

\magnification 1250
\pretolerance=500 \tolerance=1000  \brokenpenalty=5000
\mathcode`A="7041 \mathcode`B="7042 \mathcode`C="7043
\mathcode`D="7044 \mathcode`E="7045 \mathcode`F="7046
\mathcode`G="7047 \mathcode`H="7048 \mathcode`I="7049
\mathcode`J="704A \mathcode`K="704B \mathcode`L="704C
\mathcode`M="704D \mathcode`N="704E \mathcode`O="704F
\mathcode`P="7050 \mathcode`Q="7051 \mathcode`R="7052
\mathcode`S="7053 \mathcode`T="7054 \mathcode`U="7055
\mathcode`V="7056 \mathcode`W="7057 \mathcode`X="7058
\mathcode`Y="7059 \mathcode`Z="705A
\def\spacedmath#1{\def\packedmath##1${\bgroup\mathsurround =0pt##1\egroup$}
\mathsurround#1
\everymath={\packedmath}\everydisplay={\mathsurround=0pt}}
 \spacedmath{2pt}
\def\qfl#1{\buildrel {#1}\over {\longrightarrow}}

\font\eightrm=cmr8         \font\eighti=cmmi8
\font\eightsy=cmsy8        \font\eightbf=cmbx8
\font\eighttt=cmtt8        \font\eightit=cmti8
\font\eightsl=cmsl8        \font\sixrm=cmr6
\font\sixi=cmmi6           \font\sixsy=cmsy6
\font\sixbf=cmbx6\catcode`\@=11
\def\eightpoint{%
  \textfont0=\eightrm \scriptfont0=\sixrm \scriptscriptfont0=\fiverm
\def\rm{\fam\z@\eightrm}%
  \textfont1=\eighti  \scriptfont1=\sixi 
\scriptscriptfont1=\fivei
  \textfont2=\eightsy \scriptfont2=\sixsy \scriptscriptfont2=\fivesy
  \textfont\itfam=\eightit
  \def\it{\fam\itfam\eightit}%
  \textfont\slfam=\eightsl
  \def\sl{\fam\slfam\eightsl}%
  \textfont\bffam=\eightbf \scriptfont\bffam=\sixbf
  \scriptscriptfont\bffam=\fivebf
  \def\bf{\fam\bffam\eightbf}%
  \textfont\ttfam=\eighttt
  \def\tt{\fam\ttfam\eighttt}%
  \abovedisplayskip=9pt plus 3pt minus 9pt
  \belowdisplayskip=\abovedisplayskip
  \abovedisplayshortskip=0pt plus 3pt
  \belowdisplayshortskip=3pt plus 3pt 
  \smallskipamount=2pt plus 1pt minus 1pt
  \medskipamount=4pt plus 2pt minus 1pt
  \bigskipamount=9pt plus 3pt minus 3pt
  \normalbaselineskip=9pt
  \normalbaselines\rm}\catcode`\@=12
\def\note#1#2{\kern-1pt\footnote{\parindent
0.4cm$^#1$}{\vtop{\eightpoint\baselineskip12pt\hsize15.5truecm\noindent #2}}
\kern-2pt\parindent 0cm}
\def\tvi{\vrule height 2pt depth 5pt width 0pt}
\def\mono{\lhook\joinrel\mathrel{\longrightarrow}}
\def\iso{\mathrel{\mathop{\kern 0pt\smash{\longrightarrow}\tvi
}\limits^{\sim}}}
\def\tvi{\vrule height 2pt depth 5pt width 0pt}
\def\biso{\mathrel{\mathop{\kern 0pt\smash{\dasharrow}\tvi
}\limits^{\scriptscriptstyle{
\sim}}}}
\def\sdir_#1^#2{\mathrel{\mathop{\kern0pt\oplus}\limits_{#1}^{#2}}}
\def\pc#1{\tenrm#1\sevenrm}
\def\up#1{\raise 1ex\hbox{\smallf@nt#1}}
\def\tx{\kern-1.5pt -}
\def\cqfd{\kern 2truemm\unskip\penalty 500\vrule height 4pt depth 0pt width
4pt\medbreak} 
\def\virg{\raise
.4ex\hbox{,}}
\def\decale#1{\smallbreak\hskip 28pt\llap{#1}\kern 5pt}
\def\no{n\up{o}\kern 2pt}
\def\ind{\par\hskip 1truecm\relax}
\def\indp{\par\hskip 0.5truecm\relax}
\def\moins{\mathrel{\hbox{\vrule height 3pt depth -2pt width 6pt}}}
\def\rond{\kern 1pt{\scriptstyle\circ}\kern 1pt}
\def\iso{\mathrel{\mathop{\kern 0pt\longrightarrow }\limits^{\sim}}}

\def\Pic{\mathop{\rm Pic}\nolimits}

\def\rk{\mathop{\rm rk\,}\nolimits}

\def\bir{{\rm Bir}\,{\bf P}^2}
\frenchspacing
\input amssym.def
\input amssym
\vsize = 25truecm
\hsize = 16truecm
%\hoffset = -.15truecm
\voffset = -.5truecm
\parindent=0cm
\baselineskip15pt
\overfullrule=0pt

\begin
\centerline{\bf  Birational involutions of ${\bf P}^2$}
\smallskip
\smallskip \centerline{Lionel {\pc BAYLE} and Arnaud {\pc BEAUVILLE}} 
\vskip1.2cm

{\bf Introduction}\smallskip 
\ind  This paper is devoted to the classification of the elements of
order
$2$ in the group $\bir$ of birational automorphisms of ${\bf P}^2$, up to
conjugacy.  This is a classical problem, which seems to have been considered
first by Bertini [Be]. Bertini's proof is generally considered as
incomplete, as well as several other proofs which followed.  We refer to the
introduction of [C-E] for a more detailed story and for an acceptable proof.
However the result itself, as stated by these authors, is not fully
satisfactory: since they do not exclude singular fixed curves, their
classification is somewhat redundant. 

\ind We propose in this paper a  different approach, which provides a
precise and complete classification. It is based on the simple observation
that  any birational involution  of
${\bf P}^2$ is conjugate, via an appropriate birational isomorphism
$S\biso{\bf P}^2$, to a {\it
biregular} involution $\sigma$ of a rational surface
$S$. We are thus reduced to the birational classification of the pairs
$(S,\sigma)$, a problem very similar to the birational classification of
{\it real} surfaces. This classification has been achieved by classical
geometers [C], but greatly simplified in the early 80's by the introduction
of Mori theory. 
\ind In our case a direct application of Mori theory shows  that the minimal
pairs
$(S,\sigma)$ fall into two categories, those which admit a
$\sigma$\tx invariant base-point free pencil of rational curves, and those
with  $\rk
\Pic(S)^\sigma=1$. 
 The first case leads to the so-called De Jonqui\`eres involutions; in
the second case an easy lattice-theoretic argument shows that the only new
possibilities are the celebrated  Geiser and Bertini involutions. Any
birational involution is therefore conjugate to one (and only one) of these
three types.
\section{Biregular involutions of rational surfaces} 
\ind  We work over an algebraically closed field $k$ of characteristic $\not=
2$. By a {\it surface} we mean a smooth, projective, connected surface over
$k$.  

\ind  We consider pairs $(S,\sigma)$ where $S$ is a 
rational surface and $\sigma$ a non-trivial biregular involution of $S$. We
will say that $(S,\sigma)$ is {\it minimal} if any birational morphism
$f:S\rightarrow S'$ such that there exists a biregular involution
$\sigma '$ of $S'$ with $f\rond
\sigma=\sigma'\rond f$ is an isomorphism.
\ind Recall that an {\it exceptional curve} $E$ on a surface $S$ is a smooth
rational curve with $E^2=-1$.
\th Lemma
\enonce The pair $(S,\sigma)$ is minimal if and only if  each exceptional
curve $E$ on $S$ satisfies $\sigma E\not= E$ and $E\cap
\sigma E\not= \varnothing$.\endth
{\it Proof} : If $S$ contains an exceptional curve $E$ with $\sigma E= E$
(resp. $E\cap\sigma E= \varnothing$),  consider the surface $S'$ obtained
by blowing down  $E$ (resp. $E\cup\sigma E$); then
$\sigma$ induces an involution $\sigma'$ of $S'$, so that $(S,\sigma)$ is not
minimal. 
\ind Conversely, suppose that $(S,\sigma)$ is not minimal. There exists a
pair $(S',\sigma')$ and a birational morphism $f:S\rightarrow S'$ 
such that $f\rond\sigma=\sigma'\rond f$ and $f$ contracts some exceptional
curve $E$. Then $f$ contracts the divisor
$E+\sigma E$, which therefore  has negative square; this implies
$E\cdot \sigma E\le 0$, that is $\sigma E= E$ or $E\cap\sigma E=
\varnothing$.\cqfd 
\ind The only piece of Mori theory we will need is concentrated in the
following lemma; its proof follows closely that of [M], thm. 2.7.
\th Lemma
\enonce Let $(S,\sigma)$ be a minimal pair, with $\rk\Pic(S)^\sigma>1$. Then
$S$ admits a base point free pencil stable under $\sigma$. 
\endth\label{mori}
\def\rp{{\bf R}_{\scriptscriptstyle +}}
\ind Let us first recall the standard notations of Mori theory. We denote by
$NE(S)$ the  cone of effective divisor classes in $\Pic(S)\otimes{\bf R}$, by
$\overline{NE}(S)$ its closure, and by $\overline{NE}(S)_{K\ge 0}$ the
intersection of $\overline{NE}(S)$ with the half-space defined by the
condition $K_S\cdot x\ge 0$. The cone theorem ([M], 1.5 and
2.1) implies 
$$\overline{NE}(S) = \overline{NE}(S)_{K\ge 0}+  \sum_{C\in {\cal E}}
\rp[C]$$ where ${\cal E}$ is a countable set and  $C$ is a smooth
rational curve with
$C^2=-1,0$ or
$1$; moreover  if $C^2=1$  
 $S$ is  isomorphic to ${\bf P}^2$, and if $C^2=0$ $|C|$ is a base point
free pencil. 
\ind Now project the situation onto the $\sigma$\tx invariant subspace
$\Pic(S)^\sigma\otimes{\bf R}$. We get an equality (see [M], 2.6)
$$\overline{NE}(S)^\sigma = \overline{NE}(S)_{K\ge 0}^\sigma+ 
\sum_{C\in{\cal F}}
\rp[C+\sigma C]\ ,$$
where ${\cal F}$ is the subset of curves $C\in {\cal E}$ such that  the ray
$\rp[C+\sigma C]$ is extremal in $\Pic(S)^\sigma\otimes{\bf R}$. 
\ind Assume 
$\rk\Pic(S)^\sigma>1$; let $R=\rp[L]$ be an extremal ray in
$\Pic(S)^\sigma\otimes{\bf R}$. We have
$L^2\le 0$, because any element  of
$\overline{NE}(S)^\sigma$ with positive square belongs to the interior of
$NE(S)^\sigma$ ([M], Lemma 2.5). This leaves the following possibilities:
\indp $\alpha)$ $R=\rp[F]$, where $|F|$ is a base point free pencil preserved
by $\sigma $;
\indp $\beta)$ $R=\rp[E+\sigma E]$, where $E$ is an exceptional curve
and $E\cdot \sigma E=1$;
\indp $\gamma)$ $R=\rp[E+\sigma E]$, where $E$ is an exceptional curve
and $E=\sigma E$ or $E\cup \sigma E=\varnothing$. 

\ind If we  assume moreover that the pair $(S,\sigma)$ is minimal, case
$\gamma)$ does not occur. In case $\alpha)$ the conclusion is clear. In case
$\beta)$, put
$F=E+\sigma E$. We have
$F^2=0$ and
$h^0(F)\ge 2$ by Riemann-Roch; since $E$ and  $\sigma E$ do not move
linearly, this implies that $|F|$ is a base point free
pencil as required.\cqfd
\subsection\label{dp} Before stating our structure theorem for minimal pairs,
let us recall two classical examples. Let $S$ be a Del Pezzo surface of
degree 2. The linear system $|-K_S|$ defines a double covering $S\rightarrow
{\bf P}^2$, branched along a smooth quartic curve (see [D]). The involution
$\sigma$ which exchanges the two sheets of this covering is called the {\it
Geiser involution}; it satisfies $\Pic(S)^\sigma\otimes{\bf Q}\cong \Pic({\bf
P}^2)\otimes{\bf Q}={\bf Q}$. 
\ind Similarly, let $S$ be a Del Pezzo surface of degree 1. The map
$S\rightarrow {\bf P}^3$ defined by the
 linear system $|-2K_S|$ induces a degree 2 morphism of $S$ onto a
quadric cone $Q\i{\bf P}^3$, branched along the vertex $v$ of $Q$
and a smooth genus $4$ curve  [D]. The corresponding involution, the {\it
Bertini involution}, satisfies again $\rk\Pic(S)^\sigma=1$.
\th Theorem
\enonce Let $(S,\sigma)$ be a minimal pair. One of the following holds:
\indp{\rm (i)} There exists a smooth ${\bf P}^1$\tx fibration $f:S\rightarrow
{\bf P}^1$ and a non-trivial involution $\tau$ of ${\bf P}^1$ such that
$f\rond
\sigma=\tau\rond f$.
\indp{\rm (ii)} There exists a fibration  $f:S\rightarrow{\bf
P}^1$ such that $f\rond\sigma=f$; the smooth fibres of $f$ are rational
curves, on which $\sigma$ induces a non-trivial involution; any singular
fibre is the union of two rational curves exchanged by $\sigma$, meeting at
one point. 
\indp{\rm (iii)} $S$ is isomorphic to ${\bf P}^2$.
\indp{\rm (iv)} $(S,\sigma )$ is isomorphic to ${\bf P}^1\times {\bf P}^1$
with the involution
$(x,y)\mapsto (y,x)$.
\indp{\rm (v)} $S$ is a Del Pezzo surface of degree $2$ and $\sigma$ the
Geiser involution.
\indp{\rm (vi)} $S$ is a Del Pezzo surface of degree $1$ and $\sigma$ the
Bertini involution.
\endth\label{bireg}
{\it Proof} : Assume first $\rk\Pic(S)^\sigma>1$. By lemma \ref{mori} $S$
admits a $\sigma$\tx invariant pencil $|F|$ of rational curves. This
defines a fibration $f:S\rightarrow {\bf P}^1$ with fibre $F$, and an 
involution $\tau$ of ${\bf P}^1$ such that $f\rond
\sigma=\tau\rond f$.  If $f$ is smooth this gives either (i) or a particular
case of (ii). 
\ind Let $F_0$ be a singular fibre of $f$; it contains an exceptional
divisor $E$. Since $(S,\sigma)$ is minimal, we have $\sigma E\not= E$ and
$E\cdot \sigma E\ge 1$. Thus  $E+\sigma E\le F_0$ and $(E+\sigma E)^2\ge 0$,
which implies $F_0=E+\sigma E$ and $E\cdot \sigma E = 1$. 
\ind Let $p$ be the intersection point of $E$ and $\sigma E$. The involution
induced by $\sigma $ on the tangent space to $S$ at $p$ exchanges the
directions of $E$ and
$\sigma E$, hence has eigenvalues $+1$ and $-1$. It follows that there is a
fixed curve of $\sigma $ passing through $p$; this curve must be horizontal,
which forces the involution $\tau $ to be trivial. Moreover the involution
induced by $\sigma$ on  a smooth  fibre cannot be trivial, since the fixed
curve of $\sigma $ must be smooth. This gives all the properties stated in
(ii). 
\ind  Assume now
$\rk\Pic(S)^\sigma=1$. Since $\Pic(S)^\sigma$ contains an ample class, it
follows that $-K_S$ is ample, that is, $S$ is a Del Pezzo surface. If 
$\rk\Pic(S)=1$ we get case (iii). If
$\rk\Pic(S)>1$, $\sigma$ acts as
$-1$ on the orthogonal of $K_S$ in $\Pic(S)$, or in other words 
$-\sigma$ is the orthogonal reflection with respect to $K_S^\perp$. Such
a  reflection is of the form $ \displaystyle x\mapsto x- 2{(\alpha\cdot
x)\over (\alpha\cdot \alpha)}\alpha$,  with
$(\alpha\cdot \alpha)\in\{1,2\}$ and  $K_S$  proportional
to $\alpha$. If 
$K_S$ is divisible,
$S$ is isomorphic to ${\bf P}^1\times {\bf P}^1$, and we get case (iv) 
because
$\sigma $ must act non-trivially on $\Pic(S)$.
The remaining possibilities are  $K_S^2=1\ {\rm or}\ 2$. 
In these cases we have seen that the Geiser and Bertini involutions have the
required properties (\ref{dp}); they are the only ones, since an
automorphism $\gamma$ of
$S$ which acts trivially on $\Pic(S)$ is the identity (consider $S$ as
${\bf P}^2$ blown up at $9-d$ points in general position: $\gamma$
induces an automorphism of
${\bf P}^2$ which must fix these points).\cqfd
\rem{Complement}\label{ii} Let us consider case (ii) more closely. Let
$F_1,\ldots,F_s$ be the singular fibres of $f$, and $p_i\ (1\le i\le s)$ the
singular point of $F_i$. The fixed locus of $\sigma$ is a smooth curve $C$,
which passes through $p_1,\ldots,p_s$;  the covering 
$C\rightarrow {\bf P}^1$ induced by $f$ is 
of degree $2$, ramified at $p_1,\ldots,p_s$. This leads us to  distinguish
the following cases:
\indp$({\rm ii})_{sm}$: if $f$ is smooth, we have $s=0$ and $C$ is the
union of two sections of $f$ which do not intersect.
\indp$({\rm ii})_g$: if $f$ is not smooth, $C$ is a hyperelliptic curve of
genus $g\ge 0$, and $s=2g+2$.\cqfd
\ind Theorem \ref{bireg} is sufficient for our purpose, but it does not tell
us  which pairs in the  list \ref{bireg} are indeed minimal. Before 
answering that
 question  we need to work out two more
examples:\smallskip 
\rem{Examples} a) Let ${\bf F}_1$ be the surface obtained by blowing up a
point $p\in{\bf P}^2$; projecting from $p$ defines a  ${\bf P}^1$\tx
bundle $f:{\bf F}_1\rightarrow {\bf P}^1$. Any biregular involution $\sigma $
of ${\bf F}_1$ preserves this fibration, hence defines a pair $({\bf
F}_1,\sigma )$ of type (i) or $({\rm ii})_{sm}$. On the other hand
$\sigma
$ preserves the unique exceptional curve $E_1$ of ${\bf F}_1$, so the pair 
$({\bf F}_1,\sigma )$ is not minimal: $\sigma $ induces a biregular
involution of ${\bf P}^2$. In case (i) $p $ lies  on the fixed line of
$\sigma $, in case $({\rm ii})_{sm}$ it is the isolated fixed point.
\ind b) Let $Q$ be a smooth conic in ${\bf P}^2$, and $p$ a point of
${\bf P}^2\moins Q$. We define a birational involution  of ${\bf
P}^2$ by mapping a point $x$ to its harmonic conjugate on the line
$<\!p,x\!>$ with respect to the two points of $<\!p,x\!>\cap\, Q$. It is not
defined at $p$ and at the two points $q,r$ where the tangent line to $Q$
passes through
$p$. Let $S$ be the surface obtained by blowing up $p,q,r$ in ${\bf P}^2$;
the above involution extends to a biregular involution
$\sigma$ of $S$, the De Jonqui\`eres involution of degree 2. The
projection from $p$ defines a fibration $S\rightarrow {\bf P}^1$ of type
(ii) above, with 2 singular fibres. The pair
$(S,\sigma )$ is not minimal: if $E$ denotes the exceptional curve above
$p$, $\sigma E$ is the proper transform of the line $<\!q,r\!>$, so it does
not meet $E$. Blowing down $E$ and $\sigma E$ we get a pair $(T,\tau)$
with $\rk \Pic(T)=2$ and $\rk \Pic(T)^\tau=1$; inspection of the
list (\ref{bireg}) shows that it is isomorphic to ${\bf P}^1\times {\bf P}^1$
with the involution exchanging the two factors.\label{dj2}
\smallskip 
\th Proposition
\enonce The pairs $(S,\sigma )$ in the  list {\rm (\ref{bireg})} are minimal,
with the following exceptions:
\indp Case {\rm (i):} $S\cong {\bf F}_1$ 
\indp Case {\rm (ii):} $S\cong {\bf F}_1$ or  
$S$ is ${\bf P}^2$ with $3$ non-collinear points  blown up and $\sigma $ is a
De Jonqui\`eres involution of degree $2$.
\endth
{\it Proof} : The pairs $(S,\sigma )$ in (iii) to (vi) have
$\rk\Pic(S)^\sigma=1$ and therefore are minimal. The pairs $(S,\sigma )$  of
type (i) or (ii) have
$\rk\Pic(S)^\sigma=2$; thus we have to
eliminate the pairs  of these types which can be obtained 
 by blowing up either a
fixed point or two  conjugate points  in a pair $(T,\tau)$
of type (iii) to (vi).
Let $E$ be the corresponding exceptional divisor in $S$ (which may be
reducible),
$H$ the pull back to $S$ of the positive generator of
$\Pic(T)^\tau$. The group $\Pic(S)^\sigma $ is
spanned by the classes of $H$ and $E$, with $H\cdot E=0$, $H^2=1\ {\rm or}\
2$, $E^2=-1\ {\rm or}\ -2$. The
$\sigma$\tx invariant pencil $F$ is linearly equivalent to $pH-qE$ for some
integers $p,q$ which are non-negative (because
$|F|$ is  base point free) and coprime (because $F$ is not 
divisible). The condition
$F^2=0$ implies $p=q=1$, and $E^2=-H^2$. Using $F\cdot K_S=-2$ the only
possibilities are $S={\bf P}^2$ with one fixed point blown up, or $S={\bf
P}^1\times {\bf P}^1$  with the involution exchanging the factors and two
conjugate points blown up.\cqfd
\section{Birational involutions of ${\bf P}^2$}
\ind The following simple observation provides the link between biregular
involutions of rational surfaces  and birational involutions of the
plane:
\th Lemma 
\enonce Let $\iota$ be a birational involution of a surface $S_1$.
There exists a birational morphism $f:S\rightarrow S_1$ and a biregular
involution
$\sigma $ of $S$ such that $f\rond \sigma =\iota\rond f$.
\endth\label{modele}
{\it Proof} : There exists a birational morphism $f:S\rightarrow S_1$ such
that the rational map $g=\iota\rond f$ is everywhere defined (elimination of
indeterminacies, see for instance [B], II.7); moreover, $f$ is a composition
$$f:S=S_n\qfl{\varepsilon_{n-1}} S_{n-1}\longrightarrow \cdots
\longrightarrow  S_2\qfl{\varepsilon_1}S_1\ ,$$
where $ \varepsilon_i: S_{i+1}\rightarrow S_{i}$ $(1\le i\le n-1) $ is
obtained by blowing up a point $p_i\in S_{i}$. Since
$\iota $ is not defined at $p_1$, so is $g^{-1} =f^{-1} \rond \iota$; by the
universal property of blowing up [B, II.8], this implies that $g$ factors
as $S\qfl{g_1} S_2\qfl{\varepsilon_1}S_1$. Proceeding by induction we see
that $g$ factors as $f\rond \sigma $, where $\sigma $ is a birational
morphism; since $f\rond \sigma ^2=f$ $\sigma $ is an involution.\cqfd

\subsection \label{fixed}We now consider birational involutions
$\iota:S\dasharrow S$, where $S$ is a rational surface. We will
say that two  such involutions $\iota:S\dasharrow S$ and $\iota':S'\dasharrow
S'$ are birationally equivalent if there exists a birational map $\varphi :
S\dasharrow S'$ such that $\varphi \rond\iota=\iota'\rond\varphi$. 
In particular, two birational involutions of ${\bf P}^2$ are birationally
equivalent if and only if they are conjugate in the group $\bir$.
\ind Suppose
that $\iota$  fixes a curve $C$; then 
$\iota'=\varphi\rond\iota\rond\varphi^{-1} $ fixes the proper
transform of $C$ under $\varphi$, which is a curve birational to $C$ except
possibly if $C$ is rational -- in which case it may be contracted to a point.
Let us define the {\it normalized fixed curve} of $\iota$ to be the
 union of the normalizations of the non-rational curves fixed by $\iota$; it
follows from the above discussion that this is an invariant of the
birational equivalence class of $\iota$. 

\subsection\label{equi} Lemma \ref{modele} tells us that any birational
involution is birationally equivalent to a 
biregular involution $\sigma :S\rightarrow S$; moreover we can assume that
the pair
$(S,\sigma )$ is minimal. Therefore the classification of conjugacy classes
of involutions in
$\bir$ is equivalent to the classification of minimal pairs $(S,\sigma )$ up
to birational equivalence. 
We first recall the classical
examples of such involutions:\smallskip 
\rem{Examples} a) Let $S$ be a Del Pezzo surface of degree 2 and $\sigma$
the  Geiser involution (\ref{dp}). We consider $S$ as the blow up of ${\bf
P}^2$ along a set ${\cal F}$ of $7$ points in general position [D], and
denote by $\varepsilon:S\rightarrow {\bf P}^2$ the blowing up map. The
birational involution $\varepsilon\rond \sigma\rond \varepsilon^{-1} $ is
the classical Geiser involution of ${\bf P}^2$. It associates to a general
point  $x\in{\bf P}^2$ the ninth intersection point of the pencil of cubics
passing through ${\cal F}$ and $x$. The normalized fixed curve is  a 
non-hyperelliptic curve of genus
$3$.

\ind b) We define similarly the Bertini involution of ${\bf P}^2$ from the
corresponding involution on a Del Pezzo surface of degree 1 (\ref{dp}),
obtained by blowing up a set ${\cal G}$ of $8$ points in general position
in ${\bf P}^2$. It associates to a general point  $x\in{\bf P}^2$ the fixed
point of the net of sextics in ${\bf P}^2$ passing through ${\cal G}$ and
$x$ and singular along ${\cal G}$. Its
normalized fixed curve is a 
non-hyperelliptic curve of genus
 $4$, whose canonical model lies on a singular quadric.

\ind c) Let $C\i{\bf P}^2$ be a  curve of degree $d\ge 2$, and $p$ a
point of ${\bf P}^2$; we assume that $C$ has an ordinary
multiple point of multiplicity $d-2$ at $p$ and no other singularity. We
associate to $(C,p)$ the unique birational involution which preserves the
lines through $p$ and fixes the curve $C$: it maps a general point 
$x\in{\bf P}^2$ to its harmonic conjugate on the line
$<\!p,x\!>$ with respect to the two residual points of intersection of $C$
with $<\!p,x\!>$. This is a 
De Jonqui\`eres involution of degree $d$, with center
$p$ and fixed curve $C$ (the case $d=2$ was already considered in
(\ref{dj2}~b)). Its normalized fixed curve is a hyperelliptic curve\note{1}
{We consider by convention an elliptic curve as hyperelliptic.} of genus
$d-2$   for $d\ge 3$; it  is empty for
$d=2$.
\medskip 
\subsection Finally let us recall  that the ${\bf P}^1$\tx bundles over 
${\bf P}^1$ are of the form ${\bf F}_n:={\bf P}_{{\bf P}^1}({\cal O}_{{\bf
P}^1}\oplus {\cal O}_{{\bf P}^1}(n))$ for some integer $n\ge 0$. We have
${\bf F}_0={\bf P}^1\times {\bf P}^1$, and ${\bf F}_1$ is obtained by
blowing up one point in ${\bf P}^2$. For
$n\ge 1$ the fibration
$f:{\bf F}_n\rightarrow {\bf P}^1$ has a unique section $E_n$ with negative
square, and we have $E_n^2=-n$. 
\ind Let $F$ be a fibre of $f$ and $p$ a point of $F$. The {\it
elementary transformation} centered at
$p$ consists in blowing up  $p$ and blowing down the proper transform of
$F$; the surface obtained in this way is
isomorphic to ${\bf F}_{n-1}$ if  $p\notin E_n$, to ${\bf
F}_{n+1}$ if $p\in E_n$ or $n=0$.
\ind Suppose moreover that we have a birational involution $\iota$ of
${\bf F}_n$ which is regular in a neighborhood of  $F$  and
fixes
 $p$. Then after performing the elementary
transformation at $p$ we still get  a birational involution of
${\bf F}_{n\pm 1}$ which is regular in a neighborhood of the new fibre. 
\th Theorem
\enonce Every non-trivial birational involution of ${\bf P}^2$ is conjugate
to one and only one of the following:
\indp-- A De Jonqui\`eres involution of a given degree $d\ge 2$;
\indp-- A Geiser involution;
\indp-- A Bertini involution.
\endth
{\it Proof} : The unicity assertion follows from (\ref{fixed}). By
(\ref{equi}) we must prove that the involutions of the list \ref{bireg} are
birationally equivalent to one of the above types. Cases (v) and (vi) give by
definition the Geiser and Bertini involutions; 
we have seen in \ref{dj2}~b) that an involution of type (iv) is
birationally equivalent to a De Jonqui\`eres involution of degree $2$. 
\ind In case (iii), we choose a point $p\in{\bf P}^2$ with $\sigma p\not= p$
and blow up $p$ and $\sigma p$; then we blow down the proper transform of the
line $<\!p,\sigma p\!>$, which is a $\sigma $\tx invariant exceptional
curve. We obtain a pair $(T,\tau)$ with $T\cong {\bf P}^1\times {\bf P}^1$
(by stereographic projection) and $\rk\Pic(T)^\tau=1$, hence of type (iv). 
\ind In case (i), the surface $S$ is isomorphic to ${\bf F}_n$ for
some $n\ge 0$; the involution $\sigma$ has
$2$ invariant fibres, each of them containing at least $2$ fixed points. One
of these points does not lie on 
$E_n$, so performing successive elementary transformations we arrive at
$n=1$. As explained in \ref{dj2} a) we conclude that $\sigma $ is
birationally equivalent to a biregular involution of ${\bf P}^2$ (case (iv)).
\ind Let us treat case ${\rm (ii)}_{sm}$ (\ref{ii}). Again by performing 
elementary transformations we can suppose that $S$ is the surface
${\bf F}_{1}$. The fixed locus of $\sigma $ is 
the union of $E_1$ and a section which do not meet $E_1$. Blowing 
down $E_1$ gives again case (iv).
\ind It remains to  treat case ${\rm (ii)}_{g}$ for $g\ge 0$  (\ref{ii}).
Blowing down one of the components in each singular fibre we get a surface
${\bf F}_n$ with a birational involution; the fixed curve $C$ is embedded 
into ${\bf F}_n$. Performing successive elementary transformations at 
general points of $C$ leads to the same situation on 
${\bf F}_{1}$. The genus formula gives $E_1\cdot C=g$. 
 
\ind Assume that $C$ is tangent to $E_1$ at some
point $q$ of ${\bf F}_1$. Performing an elementary transformation at $q$,
then at some general point of $C$, we lower by $1$ the order of contact of
$C$ and $E_1$ at $q$. Proceeding in this way we arrive at a situation 
where $E_1$ and $C$ meet transversally at $g$ distinct points. We blow down
$E_1$ to a point
$p$ of ${\bf P}^2$; the curve $C$ maps to a plane curve $\overline{C}$ of
degree $g+2$, with an ordinary multiple point of multiplicity $g$ at $p$
and no other singularity. We get a birational involution of ${\bf P}^2$
which preserves the lines through $p$ and admits $\overline{C}$ as fixed
curve: this is the De Jonqui\`eres involution with center $p$ and fixed 
curve $\overline{C}$.\cqfd

\ind We can be more precise about the parameterization of each conjugacy
class:
\th Proposition
\enonce The map which associates to a birational involution of ${\bf P}^2$
its normalized fixed curve {\rm (\ref{fixed})} establishes a one-to-one
correspondence between:
\indp-- conjugacy classes of De Jonqui\`eres involutions
of degree $d$ and isomorphism classes of hyperelliptic curves of genus $d-2$
$(d\ge 3)$;
\indp-- conjugacy classes of Geiser involutions
 and isomorphism classes of non-hyperelliptic curves of genus $3$;
\indp-- conjugacy classes of Bertini involutions
 and isomorphism classes of non-hyperelliptic curves of genus $4$ whose
canonical model lies on a singular quadric.
\ind The De Jonqui\`eres involutions of degree $2$ form one conjugacy class.
\endth 
{\it Proof} : The result is clear for the Bertini involution: the canonical
model of the genus $3$ curve is a plane quartic; the double cover of the
plane branched along that quartic is a Del Pezzo surface of degree 2, 
which carries a canonical involution as explained in
(\ref{dp}). Similarly the  canonical model of a genus $4$ curve lies on a
unique quadric, so again we recover the Geiser involution by taking the
double cover of this quadric branched along the curve and the singular
point of the quadric.
\ind Let $g\ge 1$. A De Jonqui\`eres involution of degree $g+2$  is
determined by  a plane curve $\overline{C}$  of degree $g+2$, with an
ordinary multiple point
$p$ of multiplicity $g$ and no other singularity\note{1}{We leave to the 
reader the obvious  modifications  needed in the case $g=1$.}. 
The normalization $C$ of $\overline{C}$ is
a hyperelliptic curve of genus
$g$, with $g$ distinct points 
$p_1,\ldots ,p_g$  mapped  to $p$; the map $C\rightarrow \overline{C}\mono
{\bf P}^2$ is given by 
the linear system $|p_1+\ldots +p_g+g^1_2|$, where $g^1_2$ denotes the
degree
$2$ linear pencil  on $C$ and $p_1,\ldots ,p_g$ the points which are mapped 
to $p$. Blowing up $p$ we can view $C$ as embedded in ${\bf F}_1$; we
have $E_{1\,|C}=p_1+\ldots +p_g$. This implies in particular
$\sigma p_i\not= p_j$ for all pairs $i,j$, where  $\sigma $ stands for the
hyperelliptic involution on $C$. Let
$p_{g+1}$ be any point of $C$  such the points 
$p_1,\ldots p_{g+1},\sigma p_1,\ldots ,\sigma p_{g+1}$ are all
distinct; performing an elementary transformation at
$p_1$, then at  $\sigma p_{g+1}$, we get a  birationally
equivalent embedding $C\mono {\bf F}_1$ such that $E_{1\,|C}=p_2+\ldots
+p_{g+1}$. Continuing in this way we see that all maps of  $C$ onto a plane
curve of degree $g+2$ with one ordinary $g$\tx uple point give rise to
birationally equivalent involutions, so there is only one conjugacy class of
De Jonqui\`eres involutions with normalized fixed curve $C$. 
\ind  Finally any two degree $2$ De
Jonqui\`eres involutions are conjugate by a linear isomorphism.\cqfd

\vskip2cm
\centerline{REFERENCES} \vglue15pt\baselineskip12.8pt
\def\num#1{\smallskip\item{\hbox to\parindent{\enskip [#1]\hfill}}}
\parindent=1cm
\num{B} A. {\pc BEAUVILLE}: {\sl 	Surfaces alg\'ebriques complexes.} 
Ast\'erisque {\bf 54},  (1978).
\num{Be} E. {\pc BERTINI}: {\sl Ricerche sulle trasformazioni
univoche involutorie nel piano}. Annali di Mat. {\bf 8}, 244--286 (1877).
\num{C-E} G. {\pc CASTELNUOVO}, F. {\pc ENRIQUES}: {\sl Sulle condizioni di
razionalit\`a dei piani doppi}. Rend. del Circ. Mat. di Palermo {\bf 14},
290--302 (1900).
\num{C} A. {\pc COMESSATTI}: {\sl Fondamenti per la geometria sopra le
superficie razionali del punto de vista reale}. Math. Ann. {\bf 73}, 1--72
(1912).
\num{D} M. {\pc DEMAZURE}: {\sl Surfaces de Del Pezzo}. S\'eminaire sur les
singularit\'es des surfaces, p. 23--69; Lecture Notes {\bf 777},
Springer-Verlag (1980).
\num{M} S. {\pc MORI}: {\sl Threefolds whose canonical bundles are not
numerically effective}. Ann. of Math. {\bf 116}, 133--176 (1982).
\def\pc#1{\eightrm#1\sixrm}\eightrm\bigskip

$$\kern-0.1cm\vtop{\eightrm\hbox to 5cm{\hfill Lionel 
{\pc BAYLE}\hfill}
\hbox to 5cm{\hfill D\'epartement de Math\'ematiques\hfill}
\hbox to 5cm{\hfill UPRESA 6093 du CNRS\hfill} 
\hbox to 5cm{\hfill Universit\'e d'Angers\hfill}
\hbox to 5cm{\hfill  2 boulevard Lavoisier\hfill}
\hbox to 5cm{\hfill F-49045 {\pc ANGERS} Cedex 02\hfill}}\kern3cm
\vtop{\eightrm\hbox to 5cm{\hfill Arnaud {\pc
BEAUVILLE}\hfill}
 \hbox to 5cm{\hfill DMI -- \'Ecole Normale
Sup\'erieure\hfill} \hbox to 5cm{\hfill (UMR 8553 du CNRS)\hfill}
\hbox to 5cm{\hfill  45 rue d'Ulm\hfill}
\hbox to 5cm{\hfill F-75230 {\pc PARIS} Cedex 05\hfill}}$$
\end